\documentclass[12pt,a4]{article}%
\usepackage{amsmath}
\usepackage{graphicx}
\usepackage{hyperref}
\usepackage{eurosym}
\usepackage{amsfonts}
\usepackage{amssymb}%
\providecommand{\U}[1]{\protect\rule{.1in}{.1in}}

\begin{document}

\begin{center}
\textbf{{\large LOG-OPTIMAL AND RAPID PATHS\ IN\medskip\ } }

\textbf{{\large VON NEUMANN-GALE DYNAMICAL SYSTEMS}\textbf{{\Large \bigskip
\bigskip} } }

\textbf{\textbf{E. Babaei\footnote{Economics, University of Manchester, Oxford
Road, Manchester M13 9PL, UK. E-mail:
esmaeil.babaeikhezerloo@manchester.ac.uk.}, I.V.
Evstigneev\footnote{Economics, University of Manchester, Oxford Road,
Manchester M13 9PL, UK. E-mail: igor.evstigneev@manchester.ac.uk.
(Corresponding author.)} and K. R. Schenk-Hopp\'{e}\footnote{Economics,
University of Manchester, Oxford Road, Manchester M13 9PL, UK. E-mail:
klaus.schenk-hoppe@manchester.ac.uk.}{\Large \bigskip} } }
\end{center}

{\small \noindent\textbf{Abstract: }Von Neumann-Gale dynamical systems are
defined in terms of multivalued operators in spaces of random vectors,
possessing certain properties of convexity and homogeneity. A central role in
the theory of such systems is played by a special class of paths
(trajectories) called rapid: they grow over each time period }$t-1,t$
{\small in a sense faster than others. The paper establishes existence and
characterization theorems for such paths showing, in particular, that any
trajectory maximizing a logarithmic functional over a finite time horizon is
rapid. The proof of this result is based on the methods of convex analysis in
spaces of measurable functions. The study is motivated by the applications of
the theory of von Neumann-Gale dynamical systems to the modeling of capital
growth in financial markets with frictions -- transaction costs and portfolio
constraints.\medskip\ }

{\small \noindent\textbf{Key words and phrases:} random dynamical systems,
convex multivalued operators, von Neumann-Gale dynamical systems, rapid paths,
logarithmically optimal paths, stochastic optimization, convex analysis in
}$L_{\infty}${\small , Yosida-Hewitt decomposition, financial markets,
transaction costs, portfolio constraints, capital growth, benchmark
strategies.\medskip\ }

{\small \noindent\textbf{2010 Mathematics Subject Classifications:} 37H99,
46N10, 90C15, 91G80\medskip\ }

\textbf{\textbf{{\Large \bigskip\bigskip} {\Large \bigskip\bigskip\bigskip}
{\Large \bigskip\bigskip} {\Large \bigskip}}}

\pagebreak

\section{Introduction}

Let $(\Omega,\mathcal{F},P)$ be a complete probability space and
$\mathcal{F}_{0}\subseteq\mathcal{F}_{1}\subseteq...\subseteq\mathcal{F}%
_{N}=\mathcal{F}$ a sequence of $\sigma$-algebras containing all sets in
$\mathcal{F}$ of measure zero. For each $t=0,1,...,N$, let $X_{t}(\omega)$ be
a closed cone in an $m_{t}$-dimensional linear space $\mathbb{R}^{m_{t}}$ and
for each $t=1,...,N$, let $(\omega,a)\mapsto A_{t}(\omega,a)$ be a set-valued
operator assigning a non-empty set $A_{t}(\omega,a)\subseteq X_{t}(\omega)$ to
each $\omega\in\Omega$ and $a\in X_{t-1}(\omega)$. Throughout the paper, the
following conditions of homogeneity and convexity will be imposed on the
operator $A_{t}(\omega,\cdot)$. For each $\omega$, we have
\begin{equation}
\lambda A_{t}(\omega,a)\subseteq A_{t}\left(  \omega,\lambda a\right)
\label{r2}%
\end{equation}
for all $a\in X_{t-1}(\omega)$, $\lambda\in\lbrack0,\infty)$ and%
\begin{equation}
\theta A_{t}\left(  \omega,a\right)  +\left(  1-\theta\right)  A_{t}\left(
\omega,a^{\prime}\right)  \subseteq A_{t}\left(  \omega,\theta a+\left(
1-\theta\right)  a^{\prime}\right) \label{r3}%
\end{equation}
for all$\;a,a^{\prime}\in X_{t-1}(\omega)\ $and$\;\theta\in\left[  0,1\right]
$. (A linear combination of two sets in a vector space is the set of pairwise
linear combinations of their elements.)

The $\sigma$-algebra $\mathcal{F}_{t}$ ($t=0,...,N$) is interpreted as the
class of events occurring prior to time $t$. Vector functions of $\omega
\in\Omega$ measurable with respect to $\mathcal{F}_{t}$ represent random
vectors depending on these events. Denote for shortness by $\mathcal{L}%
_{t}^{k}$ the space $L_{\infty}\left(  \Omega,\mathcal{F}_{t},P,\mathbb{R}%
^{k}\right)  $ of essentially bounded $\mathcal{F}_{t}$-measurable functions
of $\omega\in\Omega$ with values in $\mathbb{R}^{k}$. We say that a vector
function $x(\omega)$ is a \textit{random state} of the system at time $t$ and
write $x\in\mathcal{X}_{t}$ if $x\in\mathcal{L}_{t}^{m_{t}}$ and $x(\omega)\in
X_{t}(\omega)$ almost surely (a.s.). The mappings $(\omega,a)\mapsto
A_{t}(\omega,a)$ generate a multivalued stochastic dynamical system over the
time interval $t=1,2,...,N$. A sequence of random states $x_{0}\in
\mathcal{X}_{0}$, $x_{1}\in\mathcal{X}_{1}$, $...$,$\ x_{N}\in\mathcal{X}_{N}$
is called a \textit{path }(\textit{trajectory}) of this dynamical system if
\begin{equation}
x_{t}(\omega)\in A_{t}\left(  \omega,x_{t-1}(\omega)\right)  \;\text{(a.s.)}%
.\label{r4}%
\end{equation}

Relation (\ref{r4}) can be written in the form
\begin{equation}
\left(  x_{t-1}(\omega),x_{t}(\omega)\right)  \in Z_{t}(\omega)\;\text{(a.s.)}%
,\label{r5}%
\end{equation}
where
\begin{equation}
Z_{t}(\omega)=\left\{  \left(  a,b\right)  \in X_{t-1}(\omega)\times
X_{t}(\omega):\;b\in A_{t}(\omega,a)\right\} \label{r6}%
\end{equation}
is the graph of the set-valued mapping $A_{t}\left(  \omega,\cdot\right)  $.
Clearly conditions (\ref{r2}) and (\ref{r3}) hold if and only if $Z_{t}%
(\omega)$ is a cone contained in $X_{t-1}(\omega)\times X_{t}(\omega)$. Since
$A_{t}\left(  \omega,a\right)  \neq\emptyset$ for all $a\in X_{t-1}(\omega)$,
the projection of $Z_{t}(\omega)$ on $X_{t-1}(\omega)$ coincides with
$X_{t-1}(\omega)$. It is assumed that the cones $X_{t}(\omega) $ and
$Z_{t}(\omega)$ depend $\mathcal{F}_{t}$-measurably\footnote{A set
$A(\omega)\subseteq\mathbb{R}^{k}$ is said to depend $\mathcal{F}_{t}%
$-measurably on $\omega$ if the graph $\{(\omega,a):a\in A(\omega)\}$ of the
set-valued mapping $\omega\mapsto A(\omega)$ belongs to the $\sigma$-algebra
$\mathcal{F}_{t}\times\mathcal{B(}\mathbb{R}^{k})$, where $\mathcal{B(\cdot})
$ stands for the Borel $\sigma$-algebra.} on $\omega$, which means that they
are determined by events occurring prior to time $t$.

The dynamics of the system under consideration can equivalently be described
both in terms of the mappings $A_{t}\left(  \omega,\cdot\right)  $ and in
terms of the cones $Z_{t}(\omega)$. A sequence $x_{0}\in\mathcal{X}_{0}$,
$x_{1}\in\mathcal{X}_{1}$, $...$,$\ x_{N}\in\mathcal{X}_{N}$ is a path if and
only if
\[
(x_{t-1},x_{t})\in\mathcal{Z}_{t},\ t=1,2,...,N,
\]
where%
\begin{equation}
\mathcal{Z}_{t}=\{(x,y)\in\mathcal{X}_{t-1}\times\mathcal{X}_{t}:\ \left(
x(\omega),y(\omega)\right)  \in Z_{t}(\omega)\;\text{(a.s.)}\}.\label{Z-t}%
\end{equation}

Such dynamical systems were first considered in the context of the modeling of
economic growth by von Neumann \cite{vonNeumann1937} and Gale \cite{Gale1956}.
Important contributions to the field were made by Rockafellar
\cite{Rockafellar1967}, Radner \cite{Radner1961}, McKenzie \cite{McKenzie1986}%
, Nikaido \cite{Nikaido1968}, Makarov and Rubinov \cite{MakarovRubinov1977}
and others. For reviews of this field see \cite{MakarovRubinov1977} and
\cite{EvstigneevSchenkHoppe2006}.

The classical theory of von Neumann-Gale dynamics was purely deterministic.
First attempts to build its stochastic generalization were undertaken in the
1970s by Dynkin \cite{Dynkin1971,Dynkin1972,DynkinYushkevich1979}, Radner
\cite{Radner1971} and their research groups. However, the initial attack on
the problem left many questions unanswered. Substantial progress was made only
in the late 1990s, and final solutions to the main open problems were obtained
only in the 2000s -- see \cite{EvstigneevSchenkHoppe2008}. At about the same
time it was observed \cite{DempsterEvstigneevTaksar2006} that stochastic
analogues of von Neumann-Gale dynamical systems provide a natural and
convenient framework for the modelling of financial markets with frictions
(transaction costs and portfolio constraints). This observation gave a new
momentum to studies in the field and posed new interesting problems.

In the present work, we examine the structure of paths of stochastic von
Neumann-Gale dynamical systems focusing primarily on questions of their
growth. Our main goal is to single out and investigate a class of trajectories
which grow faster in a certain sense than other trajectories over each time
period $t-1,t$. The central notion here is that of a \textit{rapid path}. Let
us give its definition. To this end we will first define the important notion
of a \textit{dual path}.

Let $X_{t}^{\ast}(\omega)$ denote the dual cone of $X_{t}(\omega)$:%
\[
X_{t}^{\ast}(\omega)=\{p\in\mathbb{R}^{m_{t}}:\ pa\geq0,\ a\in X_{t}%
(\omega)\},
\]
where $pa$ is the scalar product of the vectors $p$ and $a$ in $\mathbb{R}%
^{m_{t}}$. For shortness, we will use the notation $\mathcal{P}_{t}^{k}%
=L_{1}(\Omega,\mathcal{F}_{t},P,\mathbb{R}^{k})$ for the space of integrable
$\mathcal{F}_{t}$-measurable vector functions with values in $\mathbb{R}^{k}$.
Put%
\[
\mathcal{F}_{N+1}:=\mathcal{F}_{N}.
\]
A \textit{dual path} (\textit{dual trajectory}) is a sequence of vector
functions $p_{1}(\omega),p_{2}(\omega),...$, $p_{N+1}(\omega)$ such that
$p_{t}\in\mathcal{P}_{t}^{m_{t-1}}$ and for almost all $\omega$,
\begin{equation}
p_{t}(\omega)\in X_{t-1}^{\ast}(\omega),\ t=1,2,...,N+1,\label{p-X}%
\end{equation}
and%
\begin{equation}
\bar{p}_{t+1}(\omega)b\leq p_{t}(\omega)a\ \text{for all }(a,b)\in
Z_{t}(\omega),\ t=1,2,...,N.\label{dual-path}%
\end{equation}
Here, $\bar{p}_{t+1}(\omega):=E_{t}p_{t+1}(\omega)$ and $E_{t}(\cdot
)=E(\cdot|\mathcal{F}_{t})$ is the conditional expectation given
$\mathcal{F}_{t}$.

Let us say that a dual path $p_{1},p_{2},...,p_{N+1}$ \textit{supports} a path
$x_{0},x_{1},...,x_{N}$ if
\begin{equation}
p_{t+1}x_{t}=1,\ t=0,1,...,N\;\text{(a.s.)}.\label{p-x}%
\end{equation}
A trajectory is called \textit{rapid }if there exists a dual trajectory
supporting it. What matters in (\ref{p-x}) is that $p_{t+1}x_{t}$ is constant
(independent of time and random factors). The value $1$ for this constant is
chosen only for the sake of convenience.

The term "rapid"\ is motivated, in particular, by the fact that for each
$t=1,...,N$,%
\begin{equation}
E_{t}\frac{p_{t+1}y_{t}}{p_{t}y_{t-1}}=\frac{\bar{p}_{t+1}y_{t}}{p_{t}y_{t-1}%
}\leq\frac{\bar{p}_{t+1}x_{t}}{p_{t}x_{t-1}}=1\;\text{(a.s.)}%
\label{growth-rate}%
\end{equation}
for all paths $y_{0},y_{1},...,y_{N}$ with $p_{t}y_{t-1}>0$ (a.s.). This means
that the path $x_{0},x_{1},...,x_{N}$ maximizes the conditional expectation
given $\mathcal{F}_{t}$ of the\textit{\ growth rate }$p_{t+1}y_{t}%
/p_{t}y_{t-1}$ over each time period $(t-1,t]$, the maximum being equal to
$1$. The growth rate is measured in terms of the dual variables $p_{t}$, which
in economic and financial applications typically represent prices.

In the financial applications (on which we focus in Section 6), paths in the
dynamical system at hand represent self-financing \textit{trading strategies}.
The cones $X_{t}(\omega)$ and $Z_{t}(\omega)$ specify portfolio admissibility
constraints and self-financing constraints, respectively. Rapid paths, the
main object of our study, are counterparts of \textit{benchmark strategies}
(Platen \cite{Platen2006}, Platen and Heath \cite{PlatenHeath2006}) or
\textit{numeraire portfolios} (Long \cite{Long1990}).

This paper concentrates on the case of a finite time horizon. It generalizes
to general random cones $X_{t}(\omega)$ the results obtained in
\cite{EvstigneevFlaam1998} (also for a finite-horizon case) in a setting where
$X_{t}(\omega)$ are standard non-negative cones $\mathbb{R}_{+}^{m_{t}}$. A
central result is Theorem 1 establishing the existence of a path $x_{0}%
,x_{1},...,x_{N}$ maximizing a functional of the form $E\ln\psi(x_{N})$
(log-optimal path) and showing that this path is rapid. Extensions to an
infinite time horizon will be considered in subsequent work.

The plan of the paper is as follows. The main assumptions and results are
formulated in Section 2. In Section 3 we discuss general properties of rapid
paths. Section 4 contains some auxiliary results needed for the proof of the
main result, which is given in Section 5. Section 6 analyzes a model of a
financial market with transaction costs and portfolio constraints which is
based on von Neumann-Gale dynamical systems and to which the results of this
paper can be applied.

\section{Main results}

Let $\left\vert \cdot\right\vert $ denote the norm of a vector in a
finite-dimensional space defined as the sum of the absolute values of its
coordinates. For a finite-dimensional vector $a$, we will denote by
$\mathbb{B}(a,r)\ $the ball $\{b:|b-a|\leq r\}$. Throughout the paper it will
be assumed that conditions (\textbf{A1}) - (\textbf{A4}) we list below hold.

(\textbf{A1}) For every $t=0,1,...,N$, there exists an $\mathcal{F}_{t}%
$-measurable random vector $q_{t}(\omega)\in X_{t}^{\ast}(\omega)$ satisfying
\begin{equation}
H_{t}(\omega)^{-1}|a|\leq q_{t}(\omega)a\leq H_{t}(\omega)|a|,\ a\in
X_{t}(\omega),\ \omega\in\Omega,\label{fi}%
\end{equation}
where $H_{t}(\omega)\geq1$ is an $\mathcal{F}_{t}$-measurable function with
$E\ln H_{t}\left(  \omega\right)  <\infty$.

This condition implies, in particular, that the cone $X_{t}(\omega)$ is
\textit{pointed}, i.e., if $a\in X_{t}(\omega)$ and $-a\in X_{t}(\omega)$,
then $a=0$.

(\textbf{A2}) For every $t=1,...,N$, $\omega\in\Omega$ and $a\in
X_{t-1}(\omega)$, there exists $b\in X_{t}(\omega)$ such that $(a,b)\in
Z_{t}(\omega)$.

(\textbf{A3}) There exist constants $K_{t}$ ($t=1,...,N$) such that
$\left\vert b\right\vert \leq K_{t}\left\vert a\right\vert $ for any $(a,b)\in
Z_{t}(\omega)$ and $\omega\in\Omega$.

(\textbf{A4}) For each $t=1,2,...,N$, there exists a bounded $\mathcal{F}_{t}
$-measurable vector function $\mathring{z}_{t}=$ $(\mathring{x}_{t}%
,\mathring{y}_{t})\ $such that for all $\omega\in\Omega$, we have%
\begin{equation}
(\mathring{x}_{t}(\omega),\mathring{y}_{t}(\omega))\in Z_{t}(\omega
),\label{int-Z-y}%
\end{equation}
and%
\begin{equation}
\mathbb{B}(\mathring{y}_{t}(\omega),\varepsilon_{t})\subseteq X_{t}%
(\omega),\label{int-y}%
\end{equation}
where $\varepsilon_{t}>0$ is some constant.

For a real-valued function $\psi\left(  \omega,a\right)  $ of $\omega\in
\Omega$ and $a\in X_{t}(\omega)$ ($t=0,1,...,N$), denote by $\bar{\psi}\left(
\omega,a\right)  $ the function of $\omega\in\Omega$ and $a\in\mathbb{R}%
^{m_{t}}$ defined by%
\[
\bar{\psi}\left(  \omega,a\right)  :=\left\{
\begin{array}
[c]{cl}%
\psi\left(  \omega,a\right)  & \text{if }a\in X_{t}(\omega),\\
\infty & \text{if }a\in\mathbb{R}^{m_{t}}\setminus X_{t}(\omega),
\end{array}
\right.
\]
where "$\infty$" stands for a one-point compactification of $\mathbb{R}$.
Denote by $\Psi_{t}$ the class of real-valued functions $\psi\left(
\omega,a\right)  \geq0$ of $\omega\in\Omega$ and $a\in X_{t}(\omega)$ meeting
the following requirements:

($\psi$.1) The function $\psi\left(  \omega,\cdot\right)  $ is continuous in
$a\in X_{t}(\omega)$ for each $\omega$ and $\bar{\psi}_{t}\left(
\omega,a\right)  $ is $\mathcal{F}_{t}\times\mathcal{B}\left(  \mathbb{R}%
^{m_{t}}\right)  $-measurable in $\left(  \omega,a\right)  \in\Omega
\times\mathbb{R}^{m_{t}}$.

($\psi$.2) For all $a,a^{\prime}\in X_{t}(\omega)$, we have $\psi\left(
\omega,a+a^{\prime}\right)  $ $\geq\psi\left(  \omega,a\right)  +\psi\left(
\omega,a^{\prime}\right)  $.

($\psi$.3) The function $\psi\left(  \omega,a\right)  $ is positively
homogeneous (of degree one) in $a\in X_{t}(\omega)$:%
\[
\psi\left(  \omega,\lambda a\right)  =\lambda\psi\left(  \omega,a\right)
\text{ for any }\lambda\in\lbrack0,\infty)\text{ and }a\in X_{t}(\omega).
\]

($\psi$.4) There exists a random variable $H_{\psi}\left(  \omega\right)  >0$
such that $E\left\vert \ln H_{\psi}\left(  \omega\right)  \right\vert <\infty$
and
\begin{equation}
H_{\psi}\left(  \omega\right)  ^{-1}\left\vert a\right\vert \leq\psi\left(
\omega,a\right)  \leq H_{\psi}\left(  \omega\right)  \left\vert a\right\vert
,\;a\in X_{t}(\omega).\label{2r3}%
\end{equation}

Conditions ($\psi$.2), ($\psi$.3) and inequality (\ref{2r3}) are supposed to
hold for every $\omega\in\Omega$.

\textbf{Remark 1.} From the non-negativity of $\psi$ and requirements ($\psi
$.2), ($\psi$.3), it follows that the function $\psi\left(  \omega,a\right)
$, $a\in X_{t}(\omega)$, is concave and \textit{monotone }in$\ a$ with respect
to the partial ordering induced by the cone $X_{t}(\omega)$:%
\begin{equation}
\psi\left(  \omega,a\right)  \leq\psi\left(  \omega,a^{\prime}\right)
\;\text{if }a^{\prime}-a\in X_{t}(\omega).\label{2r3m}%
\end{equation}
Indeed, we have%
\[
\psi\left(  \omega,a^{\prime}\right)  =\psi\left(  \omega,(a^{\prime
}-a)+a\right)  \geq\psi(\omega,a^{\prime}-a)+\psi(\omega,a)\geq\psi(\omega,a).
\]

\textbf{Remark 2.} It follows from ($\psi$.4) that the expectation $E\ln
\psi\left(  \omega,x\right)  $ is well-defined and takes values in
$[-\infty,\infty)$ for any $x\in\mathcal{X}_{t}$. Furthermore, we have
\begin{equation}
E\left\vert \ln\psi\left(  \omega,x(\omega)\right)  \right\vert <\infty
\;\text{for any }x\in\text{int}\mathcal{X}_{t}.\label{2r3f}%
\end{equation}
We write $x\in$int$\mathcal{X}_{t}$ (the interior of $\mathcal{X}_{t}$) if
$\mathbb{B(}x(\omega),\varepsilon)\subseteq X_{t}(\omega)$ (a.s.) for some
constant $\varepsilon>0$. If $x\in$int$\mathcal{X}_{t}$, then $\left\vert
x(\omega)\right\vert \geq\varepsilon$ (a.s.) (since the cone $X_{t}(\omega)$
is pointed), which yields $E\ln\psi\left(  \omega,x(\omega)\right)  >-\infty$
by virtue of (\ref{2r3}). On the other hand, $E\ln\psi\left(  \omega
,x(\omega)\right)  <+\infty$, by virtue of (\ref{2r3}) and because $x(\omega)$
is essentially bounded.

\textbf{Remark 3. }Observe that the function $\psi\left(  \omega,a\right)
:=q_{t}\left(  \omega\right)  a$, where $q_{t}\left(  \omega\right)  \in
X_{t}^{\ast}(\omega)$ is the random vector described in (\textbf{A1}), belongs
to the class $\Psi_{t}$. Examples of nonlinear functions in $\Psi_{t} $ can be
constructed as follows. Let $\nu\left(  \omega,a\right)  $ be an
$\mathcal{F}_{t}\times\mathcal{B}\left(  \mathbb{R}^{m_{t}}\right)
$-measurable function of $\left(  \omega,a\right)  \in\Omega\times
\mathbb{R}^{m_{t}}$ such that $\nu\left(  \omega,\cdot\right)  $ is a norm in
$\mathbb{R}^{m_{t}}$ for each $\omega$. Since all norms in $\mathbb{R}^{m_{t}%
}$ are equivalent, there exists an $\mathcal{F}_{t}$-measurable function
$\hat{H}(\omega)\geq1$ such that
\begin{equation}
\hat{H}\left(  \omega\right)  ^{-1}\left\vert a\right\vert \leq\nu\left(
\omega,a\right)  \leq\hat{H}\left(  \omega\right)  \left\vert a\right\vert
\;\text{for all }\omega\in\Omega\text{ and }\mathbb{R}^{m_{t}}.\label{H-hat}%
\end{equation}
Assume that $E\ln\hat{H}\left(  \omega\right)  <\infty$. Define%
\begin{equation}
\psi\left(  \omega,a\right)  =q_{t}\left(  \omega\right)  a-\theta(\omega
)\nu\left(  \omega,a\right)  ,\label{psi}%
\end{equation}
where $\theta(\omega)\geq0$ is some $\mathcal{F}_{t}$-measurable function.
Clearly, the function (\ref{psi}) satisfies ($\psi$.1) - ($\psi$.3). Condition
($\psi$.4) holds if $\theta(\omega)$ is small enough. For example, take some
$\mathcal{F}_{t}$-measurable function $0<\delta(\omega)\leq1$ with
$E\left\vert \ln\delta\right\vert <\infty$. If $0\leq\theta\leq(1-\delta
)H_{t}^{-1}\hat{H}$, where $H_{t}$ is defined in (\textbf{A1}) and $\hat{H}$
in (\ref{H-hat}),\ then condition ($\psi$.4) holds with $H_{\psi}:=\delta
^{-1}H_{t}$. Indeed, we have%
\[
\delta^{-1}H_{t}\left\vert a\right\vert \geq H_{t}\left\vert a\right\vert
\geq(H_{t}-\theta\hat{H})\left\vert a\right\vert \geq q_{t}a-\theta\nu\left(
a\right)  \geq(H_{t}^{-1}-\theta\hat{H}^{-1})\left\vert a\right\vert
\]%
\[
\geq\lbrack H_{t}^{-1}-(1-\delta)H_{t}^{-1}\hat{H}\cdot\hat{H}^{-1}]\left\vert
a\right\vert =\delta H_{t}^{-1}\left\vert a\right\vert .
\]

Let $x_{0}$ be an $\mathcal{F}_{0}$-measurable vector function such that
$\mathbb{B}(x_{0}(\omega),\varepsilon_{0})\subseteq X_{0}(\omega)$, where
$\varepsilon_{0}>0$ is some constant. The random vector $x_{0}$ will be fixed
in the remainder of the paper. Denote by $\Pi\left(  x_{0},N\right)  $ the set
of paths $\xi=(x_{0},...,x_{N})$ starting from the given initial state $x_{0}
$. Fix some function $\psi_{N}(\omega,x)$ in $\Psi_{N}$ and for each path
$\xi=(x_{0},...,x_{N})\in\Pi\left(  x_{0},N\right)  $, define%
\begin{equation}
F(\xi)=E\ln\psi_{N}\left(  \omega,x_{N}(\omega)\right)  .\label{2r3a}%
\end{equation}
The main result of this paper is as follows.

\textbf{Theorem 1.} \textit{There exists a path }$\bar{\xi}=(\bar{x}_{0}%
,\bar{x}_{1},...,\bar{x}_{N})$ \textit{in }$\Pi\left(  x_{0},N\right)  $
\textit{that maximizes the functional }$F(\xi)$ \textit{over all paths }%
$\xi=(x_{0},...,x_{N})\in\Pi\left(  x_{0},N\right)  $\textit{. This path is
rapid.}

Let us say that a path is \textit{log-optimal} if it maximizes a functional of
the form (\ref{2r3a}) with some $\psi_{N}\in\Psi_{N}$. Theorem 1 shows that
any log-optimal path is rapid and thus provides an efficient method for
constructing rapid paths over a finite time horizon.

\section{General properties of rapid paths}

The first result of this section provides several equivalent definitions of a
rapid path.

Fix some $t\geq1$, $p_{t}\in\mathcal{P}_{t}^{m_{t-1}}$, $p_{t+1}\in
\mathcal{P}_{t+1}^{m_{t}}$ and $\left(  x_{t-1},x_{t}\right)  \in
\mathcal{Z}_{t}$ satisfying for almost all $\omega$
\begin{equation}
p_{t}(\omega)\in X_{t-1}^{\ast}(\omega)\text{, }p_{t+1}(\omega)\in X_{t}%
^{\ast}(\omega)\text{ and }p_{t}x_{t-1}=p_{t+1}x_{t}=1.\label{3rr1}%
\end{equation}
For any pairs $(x,y)$ of functions in $\mathcal{L}_{t}^{m_{t-1}}%
\times\mathcal{L}_{t}^{m_{t}}$ such that $(x(\omega),y(\omega))\in
Z_{t}(\omega)\ $(a.s.), consider the following four assertions.

(I) If $p_{t}x>0\ $(a.s.), then
\begin{equation}
E\left(  p_{t+1}y/p_{t}x\right)  \leq1.\label{3rr2}%
\end{equation}

(II) If $p_{t}x>0\ $(a.s.), then
\begin{equation}
E\ln\left(  p_{t+1}y/p_{t}x\right)  \leq0.\label{3rr3}%
\end{equation}

(This expectation may be a non-positive real number or $-\infty$; the function
$\ln r$ is defined as $-\infty$ for $r=0$.)

(III) The inequality%
\begin{equation}
Ep_{t+1}y\leq Ep_{t}x\label{3rr4}%
\end{equation}
holds.

(IV) With probability one, we have
\begin{equation}
E\left(  p_{t+1}(\omega)\mid\mathcal{F}_{t}\right)  b\leq p_{t}(\omega
)a\label{3rr5}%
\end{equation}
for all $(a,b)\in Z_{t}(\omega)$.

Observe that the inequalities in (\ref{3rr2}) - (\ref{3rr4}) hold as
equalities if $\left(  x,y\right)  =\left(  x_{t-1},x_{t}\right)  $.

\textbf{Proposition\ 1}. \textit{All assertions (I) - (IV) are equivalent.}

\textit{Proof.} (I)$\Rightarrow$(II). Let $(x,y)\in\mathcal{L}_{t}^{m_{t-1}%
}\times\mathcal{L}_{t}^{m_{t}}$ such that $(x(\omega),y(\omega))\in
Z_{t}(\omega)\ $(a.s.). By applying Jensen's inequality to the concave
function $\ln r$, $r\geq0$, and to the integrable non-negative random variable
$p_{t+1}y/p_{t}x$, we find $E\ln\left(  p_{t+1}y/p_{t}x\right)  \leq\ln
E\left(  p_{t+1}y/p_{t}x\right)  \leq0$, which yields (\ref{3rr3}).

(II)$\Rightarrow$(III). For any $\theta>0$, we have $\left(  x_{t-1}+\theta
x,x_{t}+\theta y\right)  \in\mathcal{L}_{t}^{m_{t-1}}\times\mathcal{L}%
_{t}^{m_{t}}$ and $\left(  x_{t-1}(\omega)+\theta x(\omega),x_{t}%
(\omega)+\theta y(\omega)\right)  \in Z_{t}(\omega)\ $(a.s.) because
$Z_{t}(\omega)$ is a convex cone. By (\ref{3rr1}) and (\ref{3rr3}),
\[
E\ln\left[  \left(  1+\theta p_{t+1}y\right)  /\left(  1+\theta p_{t}x\right)
\right]  =E\ln\left[  p_{t+1}\left(  x_{t}+\theta y\right)  /p_{t}\left(
x_{t-1}+\theta x\right)  \right]  \leq0.
\]
Consequently, $\theta^{-1}E\ln\left(  1+\theta p_{t+1}y\right)  \leq
\theta^{-1}E\ln\left(  1+\theta p_{t}x\right)  $. In the limit as
$\theta\rightarrow0$, we arrive at (\ref{3rr4}).

(III)$\Rightarrow$(IV). By virtue of the measurable selection theorem (see,
e.g., \cite{ArkinEvstigneev1987}, Appendix I), there exists a sequence
$((x_{n}(\omega),y_{n}(\omega))$, $n=1,2,...$, of measurable vector functions
such that $(x_{n},y_{n})\in\mathcal{L}_{t}^{m_{t-1}}\times\mathcal{L}%
_{t}^{m_{t}}$ and for all $\omega$ the sequence $((x_{n}(\omega),y_{n}%
(\omega))$ forms a dense subset of $Z_{t}(\omega)$. By applying (\ref{3rr4})
to $(x_{n},y_{n})$ for every $n$, we obtain $Ep_{t+1}y_{n}\leq Ep_{t}x_{n}$.
Since $Z_{t}(\omega)$ is a cone, this inequality also holds if we replace
$(x_{n},y_{n})$ by $\chi_{\Gamma}(x_{n},y_{n})$, where $\Gamma$ is any set in
$\mathcal{F}_{t}$. Then for all $(x_{n},y_{n})$, with probability one we
obtain%
\begin{equation}
E\left(  p_{t+1}\mid\mathcal{F}_{t}\right)  y_{n}\leq p_{t}x_{n}\label{E-q-p}%
\end{equation}
Let $(a,b)\in Z_{t}(\omega)$. Since the sequence $((x_{n}(\omega),y_{n}%
(\omega))$ is dense in $Z_{t}(\omega)$, it has a subsequence $(x_{n^{\prime}%
}(\omega),y_{n^{\prime}}(\omega))$ converging to $(a,b)$ and satisfying with
probability one%
\begin{equation}
E\left(  p_{t+1}\mid\mathcal{F}_{t}\right)  y_{n^{\prime}}\leq p_{t}%
x_{n^{\prime}}\label{E-q-p-n'}%
\end{equation}
for each $n^{\prime}=1,2,...$. By passing to the limit in (\ref{E-q-p-n'}) as
$n^{\prime}\rightarrow\infty$ we obtain (\ref{3rr5}) with probability one.

(IV)$\Rightarrow$(I). By applying (\ref{3rr5}) to any pairs $(x,y)$ of
functions in $\mathcal{L}_{t}^{m_{t-1}}\times\mathcal{L}_{t}^{m_{t}}$ such
that $(x(\omega),y(\omega))\in Z_{t}(\omega)\ $(a.s.) we have%

\[
E\left(  p_{t+1}\mid\mathcal{F}_{t}\right)  y\leq p_{t}x\;\text{(a.s.).}%
\]
If $p_{t}x>0\ $(a.s.), we can divide both sides of the above inequality by
$p_{t}x$ and obtain%

\[
E\left(  \left(  p_{t+1}y/p_{t}x\right)  \mid\mathcal{F}_{t}\right)
\leq1\;\text{(a.s.)}%
\]
(since $p_{t}x$ and $y$ are $\mathcal{F}_{t}$-measurable), which implies
(\ref{3rr2}).

The proof is complete.

\textbf{Proposition\ 2}.\textit{\ Replacing (\ref{dual-path}) by any of the
inequalities }$E\left(  p_{t+1}y/p_{t}x\right)  \leq1$,\textit{\ }$E\ln\left(
p_{t+1}y/p_{t}x\right)  \leq0$\textit{\ or }$Ep_{t+1}y\leq Ep_{t}%
x$\textit{\ for any pairs }$(x,y)$\textit{\ of functions in }$\mathcal{L}%
_{t}^{m_{t-1}}\times\mathcal{L}_{t}^{m_{t}}$\textit{\ such that }%
$(x(\omega),y(\omega))\in Z_{t}(\omega)\ $\textit{(a.s.), we obtain an
equivalent definition of a rapid trajectory.}

When writing the inequalities $E\left(  p_{t+1}y/p_{t}x\right)  \leq1$ and
$E\ln\left(  p_{t+1}y/p_{t}x\right)  \leq0$, we assume that the scalar product
$p_{t}x$ is strictly positive. This assumption is not needed when dealing
with\ the inequality $Ep_{t+1}y\leq Ep_{t}x$.

\textit{Proof of Proposition 2}. The assertion is a direct consequence of
Proposition 1.

\textbf{Proposition\ 3}. \textit{Let }$p_{1},p_{2},...,p_{N+1}$\textit{\ be a
dual path. For any path }$y_{0},y_{1},...$,$y_{N}$\textit{, the random
sequence }$(p_{t+1}y_{t})_{t=0}^{N}$\textit{\ is a supermartingale with
respect to the filtration} $\mathcal{F}_{1}\subseteq\mathcal{F}_{2}%
\subseteq...\subseteq\mathcal{F}_{N+1}$ \textit{and the random sequence
}$(\bar{p}_{t+1}y_{t})_{t=0}^{N}$\textit{\ is a supermartingale with respect
to the filtration }$\mathcal{F}_{0}\subseteq\mathcal{F}_{1}\subseteq
...\subseteq\mathcal{F}_{N}$.

\textit{Proof}. This is immediate from the relations:%
\[
E_{t}p_{t+1}y_{t}=\bar{p}_{t+1}y_{t}\leq p_{t}y_{t-1}\ \text{(a.s.)}%
,\ t=1,...,N,
\]
and
\[
E_{t-1}\bar{p}_{t+1}y_{t}\leq E_{t-1}p_{t}y_{t-1}=\bar{p}_{t}y_{t-1}%
\ \text{(a.s.)},\ t=1,...,N,
\]
following from (\ref{dual-path}).

\textbf{Proposition\ 4}. \textit{A path }$(x_{t})_{t=0}^{N}$\textit{\ is rapid
if and only if there exists a sequence }$(l_{t})_{t=1}^{N+1}$\textit{\ of
random vectors such that}
\begin{equation}
l_{t+1}\in X_{t}^{\ast}(\omega)\text{, }E\left\vert \ln\left(  l_{t+1}%
x_{t}\right)  \right\vert <\infty,\text{\ }l_{t+1}/l_{t+1}x_{t}\in
\mathcal{P}_{t+1}^{m_{t}}\label{3rr7}%
\end{equation}
\textit{for all }$t=0,1,...,N$,\textit{\ and}
\begin{equation}
E\ln\left(  l_{t+1}y/l_{t}x\right)  \leq E\ln\left(  l_{t+1}x_{t}/l_{t}%
x_{t-1}\right) \label{3rr8}%
\end{equation}
\textit{for any }$t=1,2,...,N$\textit{\ and for any pairs }$(x,y)$\textit{\ of
functions in }$\mathcal{L}_{t}^{m_{t-1}}\times\mathcal{L}_{t}^{m_{t}}$
\textit{such that }$(x(\omega),y(\omega))\in Z_{t}(\omega)\ $%
\textit{(a.s.)\ with }$l_{t}x>0$\textit{.}

This proposition characterizes rapid trajectories as those maximizing
\textit{the expectation of the logarithm} of the growth rate. It is important
to note that the sequence of random vectors $l_{t}$ involved in this
characterization does not necessarily satisfy the normalization condition
$l_{t+1}x_{t}=1$ (a.s.). This is in contrast with the original definition of a
rapid path, dealing with the maximization of \textit{the expectation} of the
growth rate, where the above-mentioned normalization condition is required.

\textit{Proof of Proposition 4}. If the trajectory $(x_{t})_{t=0}^{N}$ is
rapid, then we can set $l_{t}=p_{t}$, where $(p)_{t=1}^{N+1}$ is a dual path
supporting $(x_{t})_{t=0}^{N}$. The conditions contained in (\ref{3rr7}) hold
since $p_{t+1}\in X_{t}^{\ast}(\omega)$, $p_{t+1}x_{t}=1$ and $p_{t+1}%
\in\mathcal{P}_{t+1}^{m_{t}}$ for any $t=0,1,...,N$. Relation (\ref{3rr8})
turns into the inequality $E\ln\left(  p_{t+1}y/p_{t}x\right)  \leq0$, which
is true by virtue of Proposition 2.

Conversely, suppose conditions (\ref{3rr7}) and (\ref{3rr8}) are fulfilled.
Put $p_{t+1}=l_{t+1}/l_{t+1}x_{t}$\ ($t=0,1,...,N$). Note that $l_{t+1}%
x_{t}>0$ (a.s.) because $E|\ln(l_{t+1}x_{t})|<\infty$. Then for any
$t=0,1,...,N$, we have $p_{t+1}\in X_{t}^{\ast}(\omega)$, $p_{t+1}%
\in\mathcal{P}_{t+1}^{m_{t}}$ and $p_{t+1}x_{t}=1$. Since $E\left\vert
\ln\left(  l_{t+1}x_{t}/l_{t}x_{t-1}\right)  \right\vert <\infty$, it follows
from (\ref{3rr8}) that
\[
E\left[  \ln\left(  l_{t+1}y/l_{t}x\right)  -\ln\left(  l_{t+1}x_{t}%
/l_{t}x_{t-1}\right)  \right]  \leq0,
\]
which implies $E\ln\left(  p_{t+1}y/p_{t}x\right)  \leq0$. This yields
(\ref{dual-path}) by virtue of the implication (II)$\Rightarrow$(IV) proved in
Proposition 1.

The proposition is proved.

\section{Auxiliary results}

Before proving Theorem 1 we establish several lemmas.

\textbf{Lemma 1.} \textit{The functional }$F(\xi)$\textit{\ attains its
maximum over the set }$\Pi\left(  x_{0},N\right)  $ \textit{of paths }%
$(x_{0},...,x_{N})$\textit{.}

\textit{Proof.} Let us regard the class $\Pi\left(  x_{0},N\right)  $ of paths
as a subset of the space $\mathcal{L}:=\mathcal{L}_{0}^{m_{0}}\bigoplus$
$\mathcal{L}_{1}^{m_{1}}\bigoplus...$ $\bigoplus\mathcal{L}_{N}^{m_{N}}$. The
subset $\Pi\left(  x_{0},N\right)  $ is bounded in $\mathcal{L}$ by virtue of
(\textbf{A3}). Furthermore, $\Pi\left(  x_{0},N\right)  $ is convex and closed
in $\mathcal{L}$ with respect to a.s. convergence because the cones
$Z_{t}\left(  \omega\right)  $ and $X_{t}(\omega)$ are convex and closed for
each $\omega$. The functional $F\left(  \xi\right)  =E\psi_{N}\left(
\omega,x_{N}\left(  \omega\right)  \right)  $, defined for $\xi=(x_{0}%
,...,x_{N})\in\Pi\left(  x_{0},N\right)  $, is concave (which follows from the
concavity of $\psi_{N}\left(  \omega,\cdot\right)  $) and upper semicontinuous
with respect to a.s. convergence by virtue of condition ($\psi$.4) and Fatou's
lemma. This is sufficient to conclude that $F$ achieves a maximum on
$\Pi\left(  x_{0},N\right)  $; see, e.g., \cite{ArkinEvstigneev1987}, Appendix
III, Theorem 5.

The lemma is proved.

Define:%
\begin{equation}
\mathcal{U}_{t}:=\{u_{t}\in\mathcal{L}_{t}^{m_{t-1}}:\ u_{t}(\omega)\in
X_{t-1}(\omega)\text{ (a.s.)}\},\ t=1,2,...,N+1;\label{4r2a}%
\end{equation}%
\begin{equation}
\mathcal{V}_{t}:=\{v_{t}\in\mathcal{L}_{t}^{m_{t}}:\ v_{t}(\omega)\in
X_{t}(\omega)\text{ (a.s.)}\}\ [=\mathcal{X}_{t}],\ t=0,1,...,N;\label{4r2b}%
\end{equation}%
\begin{equation}
\mathcal{W}_{t}:=\{(u_{t},v_{t})\in\mathcal{U}_{t}\times\mathcal{V}%
_{t}:\left(  u_{t}\left(  \omega\right)  ,v_{t}\left(  \omega\right)  \right)
\in Z_{t}\left(  \omega\right)  \;\text{(a.s.)},\text{\ }%
t=1,...,N,\label{4r2c}%
\end{equation}
and denote by $\mathcal{W}$ the set of sequences%
\begin{equation}
\zeta=(v_{0},u_{1},v_{1},...,u_{N},v_{N},u_{N+1})\label{4r2}%
\end{equation}
with $v_{0}=x_{0}$ such that%
\begin{equation}
v_{t}\in\mathcal{V}_{t},\ t=0,...,N,\ \label{v-t}%
\end{equation}%
\begin{equation}
u_{t}\in U_{t},\ t=1,2,...,N+1,\label{u-t}%
\end{equation}%
\begin{equation}
(u_{t},v_{t})\in\mathcal{W}_{t}~,t=1,...,N,\label{w-t}%
\end{equation}
and%
\begin{equation}
E\ln\psi_{N}\left(  \omega,u_{N+1}(\omega)\right)  >-\infty.\label{4r3}%
\end{equation}
Observe that the set $\mathcal{U}_{t}$, as well as $\mathcal{X}_{t-1}$,
consists of random vectors $u_{t}(\omega)$ whose values belong to
$X_{t-1}(\omega)$ (a.s.), but these random vectors are measurable with respect
to $\mathcal{F}_{t}$ rather than $\mathcal{F}_{t-1}$, so that $\mathcal{U}%
_{t}\supset\mathcal{X}_{t-1}=\mathcal{V}_{t-1}$. Also, note that we defined
$\mathcal{F}_{N+1}$ as $\mathcal{F}_{N}$, therefore $\mathcal{U}%
_{N+1}=\mathcal{X}_{N}=\mathcal{V}_{N}$.

\textbf{Lemma 2.}\textit{\ For every }$u\in\mathcal{U}_{t}$\textit{, there
exists }$v\in\mathcal{V}_{t}$\textit{\ such that }$(u,v)\in\mathcal{W}_{t}$.

\textit{Proof.} Consider some $u\in\mathcal{U}_{t}$. Since $u(\omega)\in
X_{t-1}(\omega)$ (a.s.), we can change $u(\omega)$ on a set of measure zero
and obtain an $\mathcal{F}_{t}$-measurable vector function $u^{\prime}%
(\omega)$ such that $u^{\prime}(\omega)=u(\omega)$ (a.s.) and $u^{\prime
}(\omega)\in X_{t-1}(\omega)$ for all $\omega$. It follows from (\textbf{A2})
that for each $\omega$, there exists $b\in X_{t}(\omega)$ for which
$(u^{\prime}(\omega),b)\in Z_{t}(\omega)$. Therefore we can apply the
measurable selection theorem and construct an $\mathcal{F}_{t}$-measurable
vector $v(\omega)$ such that $(u(\omega),v(\omega))\in Z_{t}(\omega)$ (a.s.).
It follows from (\textbf{A3}) that $v(\omega)$ is essentially bounded.
Consequently, $v\in\mathcal{V}_{t}$\textit{, }$(u^{\prime},v)\in
\mathcal{W}_{t}$, and therefore $(u,v)\in\mathcal{W}_{t}$. The proof is complete.

\textbf{Lemma 3.} \textit{Let }$\zeta=(v_{0},u_{1},v_{1},...,u_{N}%
,v_{N},u_{N+1})$\textit{\ be a sequence in }$\mathcal{W}$%
\textit{\ satisfying}
\begin{equation}
v_{t-1}-u_{t}\in\mathcal{U}_{t}~,t=1,...,N+1.\label{v-U}%
\end{equation}
\textit{Then there is a path }$(y_{0},...,y_{N})$\textit{\ such that }%
$y_{0}=v_{0}$\textit{\ and }%
\[
y_{t}-v_{t}\in\mathcal{X}_{t},\ t=0,1,...,N.\mathit{\ }%
\]

\textit{Proof.} Let us proceed by induction. Put $y_{0}=v_{0}$. Suppose we
have constructed $y_{0},y_{1},...,y_{n}$ ($0\leq n<N$) satisfying%
\begin{equation}
(y_{t-1},y_{t})\in\mathcal{Z}_{t},\ 1\leq t\leq n;\label{y-Z}%
\end{equation}%
\begin{equation}
\ y_{t}-v_{t}\in\mathcal{X}_{t},\ 0\leq t\leq n.\label{y-Z-v}%
\end{equation}
(In the case of $n=0$, the constraint in (\ref{y-Z}) is absent.) Let us
construct $y_{n+1}$ for which the inclusions in (\ref{y-Z}) and (\ref{y-Z-v})
hold for $t=n+1$. Define $g_{n+1}:=y_{n}-u_{n+1}$. By virtue of (\ref{y-Z-v}),
we have $y_{n}-v_{n}\in\mathcal{X}_{n}\subseteq\mathcal{U}_{n+1}$. From
(\ref{v-U}) we get $v_{n}-u_{n+1}\in\mathcal{U}_{n+1}$. Therefore%
\[
g_{n+1}=y_{n}-u_{n+1}=(y_{n}-v_{n})+(v_{n}-u_{n+1})\in\mathcal{U}_{n+1}.
\]
By applying Lemma 2, we construct $h_{n+1}\in\mathcal{V}_{n+1}$ such that
$(g_{n+1},h_{n+1})\in\mathcal{W}_{n+1}$.

Put $y_{n+1}:=v_{n+1}+h_{n+1}$. We have $(u_{n+1},v_{n+1})\in\mathcal{W}%
_{n+1}$, $(g_{n+1},h_{n+1})\in\mathcal{W}_{n+1}$, and so%
\[
(y_{n},y_{n+1})=(u_{n+1},v_{n+1})+(g_{n+1},h_{n+1})\in\mathcal{W}_{n+1}.
\]
Since $y_{n}$ is $F_{n}$-measurable, this means that $(y_{n},y_{n+1}%
)\in\mathcal{Z}_{n+1}$, i.e., (\ref{y-Z}) holds for $t=n+1$. Furthermore%
\[
y_{n+1}-v_{n+1}=h_{n+1}\in\mathcal{V}_{n+1}=\mathcal{X}_{n+1},
\]
which gives (\ref{y-Z-v}) for $t=n+1$. Arguing by induction, we construct the
desired path $y_{0},...,y_{N}$.

The lemma is proved.

\textbf{Lemma 4.}\textit{\ There exists a sequence }$\mathring{\zeta
}=(\mathring{v}_{0},\mathring{u}_{1},\mathring{v}_{1},...,\mathring{u}%
_{N},\mathring{v}_{N},\mathring{u}_{N+1})\in\mathcal{W}$ \textit{such
that\textit{\ }}$\mathring{v}_{0}=x_{0}$\textit{, }%
\begin{equation}
\mathring{v}_{t-1}-\mathring{u}_{t}\in\text{int}\mathcal{U}_{t}%
~,t=1,...,N+1,\label{v-u-int}%
\end{equation}
\textit{and }$\mathring{u}_{N+1}\in$int$\mathcal{X}_{N}$%
.\textit{\ Furthermore, there exists a path }$\mathring{\xi}=(\mathring{x}%
_{0},\mathring{x}_{1},...,\mathring{x}_{N})\in\Pi\left(  x_{0},N\right)
$\textit{\ for which }$\mathring{x}_{t}\in$int$\mathcal{X}_{t}$, $1\leq t\leq
N$\textit{.}

Here we denote by int$~\mathcal{U}_{t}$ the interior of the set $\mathcal{U}%
_{t}$ in the topology of the space $\mathcal{L}_{t}^{m_{t-1}}$. Clearly, a
vector function $u\in\mathcal{L}_{t}^{m_{t-1}}$ belongs to int$\mathcal{U}%
_{t}$ if and only if there exists a constant $\varepsilon>0$ such that
$\mathbb{B}(u(\omega),\varepsilon)\in X_{t-1}(\omega)$ (a.s.).

\textit{Proof.} Let us argue by induction. Put $\mathring{v}_{0}=x_{0}$.
Suppose we have constructed random vectors $\mathring{v}_{t}\in\mathcal{V}%
_{t}$, $t=0,...,n$ ($0\leq n\leq N-1$) such that
\begin{equation}
\mathring{v}_{t}\in\text{int}\mathcal{U}_{t+1},\ t=0,...,n;\label{intU}%
\end{equation}
and for some $\mathring{u}_{t}\in\mathcal{U}_{t}$, $t=1,...,n$, we have%
\begin{equation}
(\mathring{u}_{t},\mathring{v}_{t})\in\mathcal{W}_{t},\ \mathring{v}%
_{t-1}-\mathring{u}_{t}\in\text{int}\mathcal{U}_{t}%
,\ t=1,...,n.\label{u-v-intU}%
\end{equation}
(For $n=0$, condition (\ref{u-v-intU}) does not make sense and is omitted.)
Let us construct $\mathring{v}_{n+1}\in\mathcal{V}_{n+1}$ and $\mathring
{u}_{n+1}\in\mathcal{U}_{n+1}\ $for which the inclusions in (\ref{intU}) and
(\ref{u-v-intU}) would hold with $t=n+1$.

Consider the pair of random vectors $(\mathring{x}_{n+1}(\omega),\mathring
{y}_{n+1}(\omega))$ described in (\textbf{A4}). Since $\mathring{x}_{n+1}%
\in\mathcal{U}_{n+1}$ and $\mathring{v}_{n}\in$int$\mathcal{U}_{n+1}$, there
exists a sufficiently small number $\lambda>0$, for which $\mathring{v}%
_{n}-\lambda\mathring{x}_{n+1}\in$int$\mathcal{U}_{n+1}$. Indeed, if
$\mathring{v}_{n}\in$int$\mathcal{U}_{n+1}$, then $\mathbb{B}(\mathring{v}%
_{n}(\omega),\delta)\subseteq X_{n}(\omega)$ (a.s.) for some $\delta>0$. By
setting $\lambda=\delta/2H$, where $H$ is a constant satisfying $\left\vert
\mathring{x}_{n+1}\right\vert \leq H$ (a.s.), we obtain that%
\[
\mathbb{B}(\mathring{v}_{n}(\omega)-\lambda\mathring{x}_{n+1}(\omega
),\delta/2)\subseteq\mathbb{B}(\mathring{v}_{n}(\omega),\delta)\subseteq
X_{n}(\omega)\ \text{(a.s.)},
\]
i.e. $\mathring{v}_{n}-\lambda\mathring{x}_{n+1}\in$int$\mathcal{U}_{n+1}$. By
defining $\mathring{v}_{n+1}:=\lambda\mathring{y}_{n+1}$ and $\mathring
{u}_{n+1}=\lambda\mathring{x}_{n+1}$, we obtain (\ref{intU}) and
(\ref{u-v-intU}) for $t=n+1$.

By applying the above induction argument, we construct a sequence
$\mathring{v}_{0},\mathring{u}_{1},\mathring{v}_{1},...,\mathring{u}%
_{N},\mathring{v}_{N}$ satisfying (\ref{intU}) and (\ref{u-v-intU}) with
$n=N$. It remains to define $\mathring{u}_{N+1}:=\mathring{v}_{N}/2$. Then
$\mathring{u}_{N+1}\in$int$\mathcal{X}_{N}$ and so $E\ln\psi_{N+1}\left(
\omega,\mathring{u}_{N+1}(\omega)\right)  >-\infty$, see Remark 2. Thus
$(\mathring{v}_{0},\mathring{u}_{1},\mathring{v}_{1},...,\mathring{u}%
_{N},\mathring{v}_{N},\mathring{u}_{N+1})\in\mathcal{W}$.

By virtue of Lemma 3 there exists a path $(\mathring{x}_{0},\mathring{x}%
_{1},...,\mathring{x}_{N})\in\Pi\left(  x_{0},N\right)  $\ for which
$\mathring{x}_{0}=x_{0}$ and $\mathring{x}_{t}-\mathring{v}_{t}\in
\mathcal{X}_{t}$,$\ t=0,1,...,N$. Since $\mathring{v}_{t-1}-\mathring{u}%
_{t}\in$int$~\mathcal{U}_{t}~$($t=1,...,N+1$) and $\mathring{u}_{t}%
\in\mathcal{U}_{t}$, we have $\mathring{v}_{t-1}\in$int$~\mathcal{U}_{t}$, and
so $\mathring{v}_{t-1}\in$int$\mathcal{X}_{t-1}$. Thus%
\[
\mathring{x}_{t}-\mathring{v}_{t}\in\mathcal{X}_{t},\ \mathring{v}_{t}%
\in\text{int}\mathcal{X}_{t},\ t=0,1,...,N,
\]
which yields $\mathring{x}_{t}=($ $\mathring{x}_{t}-\mathring{v}%
_{t})+\mathring{v}_{t}\in$int$\mathcal{X}_{t}$ for $t=0,1,...,N$.

The proof is complete.

For each sequence $\zeta=(v_{0},u_{1},v_{1},...,u_{N},v_{N},u_{N+1}%
)\in\mathcal{W}$ define
\[
G\left(  \zeta\right)  =E\ln\psi_{N}\left(  \omega,u_{N+1}(\omega)\right)
\]
and%
\[
h\left(  \zeta\right)  :=(v_{0}-u_{1},...,v_{N}-u_{N+1}).
\]
The mapping $h$ acts from the set $\mathcal{W}$ into the linear space%
\[
\mathcal{Y}:=\mathcal{L}_{1}^{m_{0}}\times\mathcal{L}_{2}^{m_{1}}%
\times...\times\mathcal{L}_{N+1}^{m_{N}}.
\]
Put%
\[
\mathcal{U}:=\mathcal{U}_{1}\times\mathcal{U}_{2}\times...\times
\mathcal{U}_{N+1}.
\]

We will show that the path $\bar{\xi}$ constructed in Lemma 1 is rapid by
analyzing the following stochastic optimization problem:

($\mathbf{P}$) Maximize the functional $G\left(  \zeta\right)  $ over the set
of sequences%
\[
\zeta=(v_{0},u_{1},v_{1},...,u_{N},v_{N},u_{N+1})\in\mathcal{W}%
\]
satisfying%
\begin{equation}
h(\zeta)\in\mathcal{U}.\label{4r4}%
\end{equation}

The following lemma shows that the path $\bar{\xi}$ generates a solution to
this problem.

\textbf{Lemma 5.} \textit{Let }$\bar{\xi}=(\bar{x}_{0},\bar{x}_{1},...,\bar
{x}_{N})$\textit{\ be a path maximizing the functional }$F(\xi)=$ $E\ln
\psi_{N}\left(  \omega,x_{N}(\omega)\right)  $\textit{\ over the set }%
$\Pi\left(  x_{0},N\right)  $ \textit{of paths }$\xi=(x_{0},x_{1},...,x_{N}%
)$\textit{. Then the sequence }$\bar{\zeta}=(\bar{x}_{0},\bar{x}_{0},\bar
{x}_{1},\bar{x}_{1},...,\bar{x}_{N-1},\bar{x}_{N},\bar{x}_{N})\in\mathcal{W}%
$\textit{\ is a solution to the optimization problem (}\textbf{P}\textit{),
and }$F(\bar{\xi})=G(\bar{\zeta})$\textit{.}

\textit{Proof.} First of all, $\bar{\zeta}\in\mathcal{W}$ since $\left(
\bar{x}_{t-1},\bar{x}_{t}\right)  \in\mathcal{Z}_{t}$ for $t=1,2,...,N$ and
\[
G(\bar{\zeta})=F(\bar{\xi})=E\ln\psi_{N}\left(  \omega,\bar{x}_{N}%
(\omega)\right)  \geq F(\mathring{\xi})=E\ln\psi_{N}(\omega,\mathring{x}%
_{N}(\omega))>-\infty
\]
(see Lemma 4 and (\ref{2r3f})). Furthermore, $h\left(  \bar{\zeta}\right)
=0\in\mathcal{U}$, so that the constraint (\ref{4r4}) is satisfied. Consider
any sequence $\zeta=(v_{0},u_{1},v_{1},...,u_{N},v_{N},u_{N+1})$ in
$\mathcal{W}$ for which $h\left(  \zeta\right)  \in\mathcal{U}$, i.e.
constraints (\ref{v-U}) hold. By virtue of Lemma 3, there is a path
$\eta=(y_{0},...,y_{N})$ such that $y_{0}=x_{0}$ and $y_{t}-v_{t}%
\in\mathcal{X}_{t}$,$\ t=0,1,...,N$. For this path, $y_{N}-v_{N}\in
\mathcal{X}_{N}$, and so%
\begin{equation}
y_{N}-u_{N+1}\in\mathcal{X}_{N}\label{y-u}%
\end{equation}
because $v_{N}-u_{N+1}\in\mathcal{X}_{N}$. Using the monotonicity of $\psi
_{N}(\omega,\cdot)$ (see (\ref{2r3m})), we obtain%
\[
G(\zeta)=E\psi_{N}(\omega,u_{N+1}(\omega))\leq E\psi_{N}(\omega,y_{N}%
(\omega))=F(\eta)\leq F(\bar{\xi})=G(\bar{\zeta}),
\]
which proves the lemma.

\section{Existence of rapid paths}

\textit{Proof of Theorem 1}. The existence of the path $\bar{\xi}=(\bar{x}%
_{0},\bar{x}_{1},...,\bar{x}_{N})$ maximizing the functional (\ref{2r3a}) was
established in Lemma 1. To show that $\bar{\xi}$ is rapid we will apply to the
optimization problem (\textbf{P})\ a general version of the Kuhn-Tucker
theorem established, e.g., in \cite{Hurwicz1958}, Theorem 5.3.1. The set
$\mathcal{W}$ is convex. The set $\mathcal{U}$ is a convex cone. The mapping
$g$ is linear. The functional $G\left(  \zeta\right)  $, $\zeta\in\mathcal{W}%
$, is concave and takes on real values. Thus, in order to justify the use of
the Kuhn-Tucker theorem we have to check Slater's condition:

(\textbf{S}) There is an element $\mathring{\zeta}\in\mathcal{W}$ such that
$g(\mathring{\zeta})$ belongs to the interior int$\mathcal{U}$ of the cone
$\mathcal{U}$ in the topology of the space $\mathcal{Y}$.

This condition holds because the sequence $\mathring{\zeta}\in\mathcal{W}$
constructed in Lemma 4 possesses the properties listed in (\textbf{S}):
relations (\ref{u-v-intU}) mean that $g(\mathring{\zeta})\in$int$\mathcal{U}$.

By the Kuhn-Tucker theorem applied to problem (\textbf{P}), there exists a
continuous linear functional $\pi$ on the space $\mathcal{Y}$ such that
$\left\langle \pi,y\right\rangle \geq0$ for $y\in\mathcal{U}$ and
\[
G\left(  \zeta\right)  +\left\langle \pi,h\left(  \zeta\right)  \right\rangle
\leq G\left(  \bar{\zeta}\right)
\]
for any $\zeta\in\mathcal{W}$. The functional $\pi$ can be represented in the
form $\pi=\left(  \pi_{1},...,\pi_{N+1}\right)  $, where $\pi_{t}$ is a
continuous linear functional on the space $\mathcal{L}_{t}^{m_{t-1}}%
=L_{\infty}(\Omega,\mathcal{F}_{t},P,\mathbb{R}^{m_{t-1}})$ such that $\pi
_{t}\in\left(  \mathcal{U}_{t}\right)  ^{\ast}$, i.e.,
\begin{equation}
\left\langle \pi_{t},u_{t}\right\rangle \geq0\ \text{for }u_{t}\in
\mathcal{U}_{t}\ (t=1,...,N+1).\label{non-neg}%
\end{equation}
Thus we have%
\[
E\ln\psi_{N}\left(  \omega,u_{N+1}(\omega)\right)  +\sum_{t=1}^{N+1}%
\left\langle \pi_{t},v_{t-1}-u_{t}\right\rangle
\]%
\begin{equation}
\leq E\ln\psi_{N}\left(  \omega,\bar{x}_{N}(\omega)\right)  \ [=G(\bar{\zeta
})=F(\bar{\xi})].\label{4r5}%
\end{equation}
for any $\zeta=(v_{0},u_{1},v_{1},...,u_{N},v_{N},u_{N+1})\in\mathcal{W}$.

By virtue of the Yosida-Hewitt theorem \cite{YosidaHewitt1952}, each of the
functionals $\pi_{t}$ can be decomposed into the sum $\pi_{t}=$ $\pi_{t}^{a}+
$ $\pi_{t}^{s}$ of two functionals $\pi_{t}^{a}$, $\pi_{t}^{s}\in$ $\left(
\mathcal{L}_{t}^{m_{t-1}}\right)  ^{\ast}$, where $\pi_{t}^{a}$ is
\textit{absolutely continuous} and $\pi_{t}^{s}$ is \textit{singular}.
According to the definitions of $\pi_{t}^{a}$ and $\pi_{t}^{s}$, there is a
vector function $p_{t}\in\mathcal{P}_{t}^{m_{t-1}}=L_{1}(\Omega,\mathcal{F}%
_{t},P,\mathbb{R}^{m_{t-1}})$ such that%
\begin{equation}
\left\langle \pi_{t}^{a},y\right\rangle =Ep_{t}y,\ y\in\mathcal{L}%
_{t}^{m_{t-1}},\label{Abs-c}%
\end{equation}
and there exist sets $\Gamma_{t}^{1}\supseteq\Gamma_{t}^{2}\supseteq...$ in
$\mathcal{F}_{t}$ for which $P\left(  \Gamma_{t}^{k}\right)  \rightarrow0$ as
$k\rightarrow\infty$ and
\begin{equation}
\left\langle \pi_{t}^{s},y\chi_{\Gamma_{t}^{k}}\right\rangle =\left\langle
\pi_{t}^{s},y\right\rangle ,\ y\in\mathcal{L}_{t}^{m_{t-1}},\label{Sing}%
\end{equation}
where $\chi_{\Gamma_{t}^{k}}$ is the indicator function of the set $\Gamma
_{t}^{k}$.

Observe that relation (\ref{non-neg}) remains valid if we replace $\pi_{t}$ by
$\pi_{t}^{a}$. Indeed, the inclusion $u_{t}\in\mathcal{U}_{t}$ means that
$u_{t}\in\mathcal{L}_{t}^{m_{t-1}}$ and$\ u_{t}(\omega)\in X_{t-1}(\omega)$
(a.s.). This implies $u_{t}\chi_{\Delta_{t}^{k}}\in\mathcal{L}_{t}^{m_{t-1}}$
and $u_{t}(\omega)\chi_{\Delta_{t}^{k}}(\omega)\in X_{t-1}(\omega)$ (a.s.),
where $\Delta_{t}^{k}:=$ $\Omega\setminus\Gamma_{t}^{k}$. Consequently, we
have%
\begin{equation}
0\leq\left\langle \pi_{t},u_{t}\chi_{\Delta_{t}^{k}}\right\rangle
=\left\langle \pi_{t}^{a},u_{t}\chi_{\Delta_{t}^{k}}\right\rangle
+\left\langle \pi_{t}^{s},u_{t}\chi_{\Delta_{t}^{k}}\right\rangle
=\left\langle \pi_{t}^{a},u_{t}\chi_{\Delta_{t}^{k}}\right\rangle
,\label{non-neg1}%
\end{equation}
where $\left\langle \pi_{t}^{a},u_{t}\chi_{\Delta_{t}^{k}}\right\rangle
\rightarrow\left\langle \pi_{t}^{a},u_{t}\right\rangle $. By passing to the
limit in (\ref{non-neg1}), we obtain that%
\begin{equation}
Ep_{t}u_{t}=\left\langle \pi_{t}^{a},u_{t}\right\rangle \geq0,\ u_{t}%
\in\mathcal{U}_{t}.\label{pi-a-U}%
\end{equation}
Since $\mathcal{U}_{t}$ consists of those functions $u_{t}\in\mathcal{L}%
_{t}^{m_{t-1}}$ for which $u_{t}(\omega)\in X_{t-1}(\omega)$ (a.s.),
inequality (\ref{pi-a-U}) implies $p_{t}(\omega)\in X_{t-1}^{\ast}(\omega
)\ $(a.s.)$,\ t=1,2,...,N+1$, i.e. condition (\ref{p-X}) holds.

Furthermore, if (\ref{non-neg}) holds, then
\begin{equation}
\left\langle \pi_{t}^{s},u_{t}\right\rangle \geq0\ \text{for }u_{t}%
\in\mathcal{U}_{t}\ (t=1,...,N+1).\label{non-neg-s}%
\end{equation}
Indeed, by virtue of (\ref{Sing}) we have%
\[
0\leq\left\langle \pi_{t},u_{t}\chi_{\Gamma_{t}^{k}}\right\rangle
=\left\langle \pi_{t}^{a},u_{t}\chi_{\Gamma_{t}^{k}}\right\rangle
+\left\langle \pi_{t}^{s},u_{t}\chi_{\Gamma_{t}^{k}}\right\rangle
=\left\langle \pi_{t}^{a},u_{t}\chi_{\Gamma_{t}^{k}}\right\rangle
+\left\langle \pi_{t}^{s},u_{t}\right\rangle ,
\]
and so%
\[
\left\langle \pi_{t}^{s},u_{t}\right\rangle \geq-\left\langle \pi_{t}%
^{a},u_{t}\chi_{\Gamma_{t}^{k}}\right\rangle \rightarrow0,
\]
which proves (\ref{non-neg-s}).

Let us show that relation (\ref{4r5}) remains valid if we replace $\pi_{t}$ by
$\pi_{t}^{a}$. We will prove this by way of induction. Fix some $\zeta
=(v_{0},u_{1},v_{1},...,u_{N},v_{N},u_{N+1})\in\mathcal{W}$ and consider the
inequality%
\[
E\ln\psi_{N}\left(  \omega,u_{N+1}(\omega)\right)  +\sum_{t=1}^{N+1}%
\left\langle \pi_{t}^{a},v_{t-1}-u_{t}\right\rangle
\]%
\begin{equation}
+\sum_{t=1}^{M}\left\langle \pi_{t}^{s},v_{t-1}-u_{t}\right\rangle \leq
E\ln\psi_{N}\left(  \omega,\bar{x}_{N}(\omega)\right) \label{4r6}%
\end{equation}
where $M\in\left\{  0,...,N+1\right\}  $. If $M=0$, then the second sum in
(\ref{4r6})\ is formally defined as 0. For $M=N+1$, relation (\ref{4r6}) is
equivalent to (\ref{4r5}). Suppose inequality (\ref{4r6}) is true for some
$M\in$ $\left\{  1,...,N+1\right\}  $. Let us show that this inequality is
true for $M-1$.

Since $P\left(  \Gamma_{M}^{k}\right)  \rightarrow0$ as $k\rightarrow\infty$,
we can find a sequence of real numbers $\epsilon_{k}\in\left(  0,1\right)  $
such that $\epsilon_{k}\rightarrow0$ and
\begin{equation}
P\left(  \Gamma_{M}^{k}\right)  \ln\epsilon_{k}\rightarrow0\label{gamma}%
\end{equation}
(e.g., we can define $\epsilon_{k}:=\exp(-\mu_{k}^{-1/2})$). Put
\[
\Delta_{M}^{k}=:\Omega\setminus\Gamma_{M}^{k};\ \gamma_{M}^{k}\left(
\omega\right)  :=\epsilon_{k}\chi_{\Gamma_{M}^{k}}\left(  \omega\right)
+\chi_{\Delta_{M}^{k}}\left(  \omega\right)  ;
\]%
\[
(u_{t}^{k},v_{t}^{k}):=\gamma_{M}^{k}(u_{t},v_{t}),\ M\leq t\leq
N,\ u_{N+1}:=\gamma_{M}^{k}u_{N+1};
\]%
\[
v_{0}^{k}=:v_{0},\;(u_{t}^{k},v_{t}^{k}):=(u_{t},v_{t}),\ 1\leq t<M;\;
\]
and
\[
\zeta^{k}:=(v_{0}^{k},u_{1}^{k},v_{1}^{k},...,u_{N}^{k},v_{N}^{k},u_{N+1}%
^{k})\ (k=1,2,...).
\]
We can see that
\[
\left(  u_{t}^{k}\left(  \omega\right)  ,v_{t}^{k}\left(  \omega\right)
\right)  \in Z_{t}\left(  \omega\right)  \text{\ (a.s.)},
\]
and so $(u_{t}^{k},v_{t}^{k})\in\mathcal{W}_{t}$. Furthermore,%
\[
v_{0}^{k}=v_{0}=x_{0},\ u_{N+1}^{k}=\gamma_{M}^{k}u_{N+1}\in\mathcal{U}%
_{N+1}=\mathcal{X}_{N}=\mathcal{V}_{N},
\]
and%
\[
E\ln\psi_{N}\left(  \omega,u_{N+1}^{k}(\omega)\right)  =E\ln\gamma_{M}%
^{k}(\omega)+E\ln\psi_{N}\left(  \omega,u_{N+1}(\omega)\right)
\]%
\begin{equation}
=P\left(  \Gamma_{M}^{k}\right)  \ln\epsilon_{k}+E\ln\psi_{N}\left(
\omega,u_{N+1}(\omega)\right)  ,\label{4r7}%
\end{equation}
therefore $E\ln\psi_{N}\left(  \omega,u_{N+1}^{k}(\omega)\right)  >-\infty$.
Consequently, the sequence $\zeta^{k}$ belongs to $\mathcal{W}$, and we can
apply inequality (\ref{4r6}) to this sequence.

This yields%
\[
E\ln\psi_{N}\left(  \omega,u_{N+1}^{k}(\omega)\right)  +\sum_{t=1}%
^{N+1}\left\langle \pi_{t}^{a},v_{t-1}^{k}-u_{t}^{k}\right\rangle +\sum
_{t=0}^{M-1}\left\langle \pi_{t}^{s},v_{t-1}-u_{t}\right\rangle +
\]%
\begin{equation}
\left\langle \pi_{M}^{s},v_{M-1}-u_{M}^{k}\right\rangle \leq E\ln\psi
_{N}\left(  \omega,\bar{x}_{N}(\omega)\right)  .\label{4r8}%
\end{equation}
Here, the vectors $u_{t}^{k}$, $v_{t}^{k}$ are uniformly bounded, and we have
$u_{t}^{k}\rightarrow u_{t}$, $v_{t}^{k}\rightarrow v_{t}$ (a.s.).
Consequently,
\[
\left\langle \pi_{t}^{a},v_{t-1}^{k}-u_{t}^{k}\right\rangle =Ep_{t}\left(
v_{t-1}^{k}-u_{t}^{k}\right)  \rightarrow Ep_{t}\left(  v_{t-1}-u_{t}\right)
\]
for each $t$. By virtue of (\ref{4r7}) and (\ref{gamma}), $E\ln\psi_{N}\left(
\omega,u_{N}^{k}(\omega)\right)  \rightarrow$ $E\ln\psi_{N}\left(
\omega,u_{N}(\omega)\right)  $. Furthermore,%
\[
\left\langle \pi_{M}^{s},v_{M-1}-u_{M}^{k}\right\rangle \geq-\left\langle
\pi_{M}^{s},u_{M}^{k}\right\rangle =-\left\langle \pi_{M}^{s},\chi_{\Gamma
_{M}^{k}}u_{M}^{k}\right\rangle =-\epsilon_{k}\left\langle \pi_{M}^{s}%
,u_{M}\right\rangle \rightarrow0,
\]
where the inequality in this chain of relations follows from (\ref{non-neg-s}%
)\ because $v_{M-1}\in V_{M-1}\subseteq U_{M}$. Thus, by passing to the limit
in (\ref{4r8}), we conclude that inequality (\ref{4r6}) remains true if we
replace $M$ by $M-1$.

We have constructed a sequence of functions $p_{t}\in\mathcal{P}_{t}$,
$t=1,...,N+1$, such that
\[
E\ln\psi_{N}\left(  \omega,u_{N+1}(\omega)\right)  +\sum_{t=1}^{N+1}%
Ep_{t}(v_{t-1}-u_{t})
\]%
\begin{equation}
\leq E\ln\psi_{N}\left(  \omega,\bar{x}_{N}(\omega)\right)  \ [=G(\bar{\zeta
})=F(\bar{\xi})]\label{ineq-p}%
\end{equation}
for any $\zeta=(v_{0},u_{1},v_{1},...,u_{N},v_{N},u_{N+1})\in\mathcal{W}$. Let
us show that $(p_{1},...,p_{N+1})$ is a dual path supporting the path
$\bar{\xi}$. We can write the sum in (\ref{ineq-p}) as%
\[
\sum_{t=1}^{N+1}Ep_{t}(v_{t-1}-u_{t})=-Ep_{N+1}u_{N+1}+\sum_{t=1}^{N}E\left(
p_{t+1}v_{t}-p_{t}u_{t}\right)  +Ep_{1}v_{0},
\]
where $Ep_{1}v_{0}=Ep_{1}x_{0}=Ep_{1}\bar{x}_{0}$, and since%
\[
-Ep_{N+1}\bar{x}_{N}+\sum_{t=1}^{N}E\left(  p_{t+1}\bar{x}_{t}-p_{t}\bar
{x}_{t-1}\right)  +Ep_{1}\bar{x}_{0}=0,
\]
we can see that inequality (\ref{ineq-p}) implies
\[
E\ln\psi_{N}\left(  \omega,u_{N+1}(\omega)\right)  -Ep_{N+1}u_{N+1}+\sum
_{t=1}^{N}E\left(  p_{t}v_{t-1}-p_{t}u_{t}\right)  +Ep_{1}x_{0}%
\]%
\[
\leq E\ln\psi_{N}\left(  \omega,\bar{x}_{N}(\omega)\right)  -Ep_{N+1}\bar
{x}_{N}+\sum_{t=1}^{N}E\left(  p_{t+1}\bar{x}_{t}-p_{t}\bar{x}_{t-1}\right)
+Ep_{1}x_{0}.
\]
In turn, this yields%
\begin{equation}
E\left(  p_{t+1}v_{t}-p_{t}u_{t}\right)  \leq E\left(  p_{t+1}\bar{x}%
_{t}-p_{t}\bar{x}_{t-1}\right)  ,\ (u_{t},v_{t})\in\mathcal{W}_{t}%
,\ t=1,...,N,\label{4r-p2}%
\end{equation}
and%
\begin{equation}
E\ln\psi_{N}\left(  u\right)  -Ep_{N+1}u\leq E\ln\psi_{N}\left(  \bar{x}%
_{N}\right)  -Ep_{N+1}\bar{x}_{N},\ u\in\mathcal{U}_{N+1}=\mathcal{X}%
_{N}.\label{4r-p1}%
\end{equation}

Recall that $\mathcal{W}_{t}$ consists of $(u_{t},v_{t})\in\mathcal{L}%
_{t}^{m_{t-1}}\times\mathcal{L}_{t}^{m_{t}}$ for which $(u_{t}(\omega
),v_{t}(\omega))\in Z_{t}(\omega)$ (a.s.). The set $Z_{t}\left(
\omega\right)  $ is a cone, consequently, inequality (\ref{4r-p2}) will remain
valid if we multiply $\left(  u_{t},v_{t}\right)  $ by any positive constant.
This and the fact that $(\bar{x}_{t-1}\left(  \omega\right)  ,\bar{x}%
_{t}\left(  \omega\right)  )\in Z_{t}\left(  \omega\right)  $ (a.s.), yields%
\begin{equation}
Ep_{t+1}v_{t}-Ep_{t}u_{t}\leq0=Ep_{t+1}\bar{x}_{t}-Ep_{t}\bar{x}%
_{t-1},\ t=1,...,N,\label{4r-p4}%
\end{equation}
for all pairs $(u_{t},v_{t})$ of functions in $\mathcal{L}_{t}^{m_{t-1}}%
\times\mathcal{L}_{t}^{m_{t}}$ such that $\left(  u_{t}\left(  \omega\right)
,v_{t}\left(  \omega\right)  \right)  \in Z_{t}\left(  \omega\right)  $
(a.s.). By virtue of Proposition 2, this means that $(p_{1},...,p_{N+1})$ is a
dual path.

It remains to prove that the dual path we have constructed supports the path
$\bar{\xi}$. To this end let us show, by using (\ref{4r-p1}), that
$p_{N+1}\bar{x}_{N}=1$ (a.s.). Put
\[
\gamma:=p_{N+1}\bar{x}_{N},\ \delta_{k}:=[\gamma+k^{-1}]^{-1},\ u^{k}%
:=\delta_{k}\bar{x}_{N}.
\]
The random variable $\gamma$ is a.s. non-negative, $\mathcal{F}_{N}%
$-measurable (because $\mathcal{F}_{N}=\mathcal{F}_{N+1}$) and integrable
because $p_{N+1}$ is integrable and $\bar{x}_{N}$ is essentially bounded. The
random variable $\delta_{k}\ $is $\mathcal{F}_{N}$-measurable and satisfies
$0\leq\delta_{k}\leq k$ (a.s.). Consequently, $\ u^{k}(\omega)\in X_{N}%
(\omega)$ (a.s.), and we can apply inequality (\ref{4r-p1}) to $u=u^{k}$. We
have%
\[
E\ln\psi_{N}\left(  u^{k}\right)  -Ep_{N+1}u^{k}=E\ln\delta_{k}+E\ln\psi
_{N}\left(  \bar{x}_{N}\right)  -E\delta_{k}\gamma.
\]
Since $-\infty<E\ln\psi_{N}\left(  \bar{x}_{N}\right)  <+\infty$, it follows
from this equality and (\ref{4r-p1}) that%
\[
E\ln\delta_{k}\leq E\delta_{k}\gamma-E\gamma.
\]
We have $0\leq\delta_{k}\gamma=\gamma\lbrack\gamma+k^{-1}]^{-1}\leq1$ (a.s.),
and so $\lim E\delta_{k}\gamma=1$. By using Fatou's lemma, we obtain%
\[
-E\ln\gamma=E\lim\inf\ln\delta_{k}\leq\lim\inf E\ln\delta_{k}\leq1-E\gamma.
\]
The use of Fatou's lemma is justified because
\[
\ln\delta_{k}=-\ln(\gamma+k^{-1})\geq-\ln(\gamma+1)\geq-\gamma,
\]
where $\gamma$ is integrable. Thus we obtain $-E\ln\gamma\leq1-E\gamma$, or
equivalently, $E(\gamma-1-\ln\gamma)\leq0$. At the same time, we always have
$\gamma-1-\ln\gamma\geq0$. Therefore $p_{N+1}\bar{x}_{N}=\gamma=1$ (a.s.).

By virtue of the equality in (\ref{4r-p4}), we have%
\begin{equation}
Ep_{1}\bar{x}_{0}=Ep_{2}\bar{x}_{1}=...=Ep_{N+1}\bar{x}_{N}.\label{4r-p5}%
\end{equation}
We can replace in (\ref{4r-p4}) $\left(  u_{t},v_{t}\right)  $ by
$\chi_{\Gamma}\left(  u_{t},v_{t}\right)  $, where $\Gamma$ is any set in
$\mathcal{F}_{t}$. This yields
\begin{equation}
E(p_{t+1}v_{t}~|\mathcal{F}_{t})-p_{t}u_{t}\leq0\text{ (a.s.)}\label{4r-p6}%
\end{equation}
for all $(u_{t},v_{t})\in\mathcal{L}_{t}^{m_{t-1}}\times\mathcal{L}_{t}%
^{m_{t}}$ such that $\left(  u_{t}\left(  \omega\right)  ,v_{t}\left(
\omega\right)  \right)  \in Z_{t}\left(  \omega\right)  $ (a.s.). By using
(\ref{4r-p6}) and (\ref{4r-p5}), we obtain $E\left(  p_{t+1}\bar{x}%
_{t}|\mathcal{F}_{t}\right)  =p_{t}\bar{x}_{t-1}$, $t=1,...,N$. Since
$p_{N+1}\bar{x}_{N}=1$, we conclude that $p_{N+1}\bar{x}_{N}=$ $p_{N}\bar
{x}_{N-1}=...=$ $p_{1}\bar{x}_{0}=1$.

Theorem 1 is proved.

In the course of the above proof, we used a procedure for deriving necessary
conditions for an extremum based on the Yosida-Hewitt theorem. Apparently the
first who applied such methods in optimization (in the context of
continuous-time optimal control) were Dubovitskii and Milyutin
\cite{DubovitskiiMilyutin1968}. Techniques based on the Yosida-Hewitt theorem
were used in the analysis of discrete-time stochastic models of economic
dynamics and related problems of stochastic programming by Radner
\cite{Radner1972,Radner1973}, Evstigneev \cite{Evstigneev1972,Evstigneev1976}%
,\ Rockafellar and Wets\ \ \cite{RockafellarWets1975}, and others. For a
comprehensive review of early literature in this field see the book by Arkin
and Evstigneev \cite{ArkinEvstigneev1987}.

\section{A financial market model}

In this section we give an example of a model for a financial market with
transaction costs and portfolio constraints that can be included in the
framework of von Neumann-Gale dynamical systems. We provide conditions that
guarantee the validity of assumptions (\textbf{A1})-(\textbf{A4}) introduced
above, which makes it possible to apply in this context the results of the
present paper. As regards the aspect of growth-optimal investments, the model
under consideration extends the one studied in Bahsoun et al.
\cite{Bahsoun-et-al2013}. In the latter, portfolio constraints are specified
by the cones $X_{t}=\mathbb{R}_{+}^{m}$, i.e., short-selling is not allowed.
Here, short sales are permitted, but are subject to certain
constraints---margin requirements (see below). In this paper, we only briefly
discuss the financial aspects, referring the reader for details to
\cite{Bahsoun-et-al2013}.

We consider a market where $m$ assets are traded at dates $t=1,2,\ldots,N$.
Random vectors $a(\omega)\in\mathbb{R}^{m}$ are interpreted as (contingent)
portfolios of assets. Positions $a_{i}(\omega)$ of the portfolio
$a(\omega)=(a_{1}(\omega),...,a_{m}(\omega))\in\mathbb{R}^{m}$ are measured in
terms of their values in the market prices. For each $t=0,1,\ldots,N$ and
$i=1,\ldots,m$ the following $\mathcal{F}_{t}$-measurable random variables are
given: asset prices $S_{t,i}(\omega)>0$ and transaction cost rates for selling
and buying assets $0\leq\lambda_{t,i}^{+}(\omega)<1$, $\lambda_{t,i}%
^{-}(\omega)\geq0$. We denote by $R_{t,i}=S_{t,i}/S_{t-1,i}$ the (gross)
return on asset $i$ at time $t$. We omit $\omega$ in the notation where it
does not lead to ambiguity.

The portfolio constraints in the model are specified by the cones
\begin{equation}
X_{t}(\omega)=\left\{  a\in\mathbb{R}^{m}:\sum_{i=1}^{m}(1-\lambda_{t,i}%
^{+}(\omega))a_{+}^{i}\geq\mu_{t}\sum_{i=1}^{m}(1+\lambda_{t,i}^{-}%
(\omega))a_{-}^{i}\right\}  ,\label{Xt-example}%
\end{equation}
where $\mu_{t}>1$ are constants (independent of $\omega$). The inequalities in
(\ref{Xt-example}) express \textit{margin requirements}: for an admissible
portfolio, the total value of its long positions must cover a \textit{margin}
$\mu_{t}$ (in the U.S. equity markets $\mu_{t}=1.5$) times the total value its
short positions. These values are computed taking into account transaction
costs for buying and selling assets.

Trading in the market proceeds as follows. At each date $t=1,2,\ldots,N$ a
trader can rebalance her portfolio $a(\omega)\in X_{t-1}(\omega)$ purchased at
the previous date $t-1$ to a new portfolio $b(\omega)\in X_{t}(\omega)$. The
possibilities of rebalancing are specified by the inequality $\psi_{t}%
(\omega,a,b)\geq0$, where
\[
\psi_{t}(a,b)=\sum_{i=1}^{m}(1-\lambda_{t,i}^{+})(R_{t,i}a^{i}-b^{i})_{+}%
-\sum_{i=1}^{m}(1+\lambda_{t,i}^{-})(R_{t,i}a^{i}-b^{i})_{-}.
\]
The first sum represents the amount of money the trader receives for selling
assets, the second sum is the amount of money she pays for buying assets,
including transaction costs. The inequality $\psi_{t}(a,b)\geq0$ means that
the trader does not use external funds to rearrange her portfolio, and so it
can be regarded as a \textit{self-financing condition}.

Define
\begin{equation}
{Z}_{t}(\omega)=\left\{  (a,b)\in X_{t-1}(\omega)\times X_{t}(\omega):\psi
_{t}(\omega,a,b)\geq0\right\}  .\label{Example2-Gt1}%
\end{equation}
Observe that ${Z}_{t}(\omega)$ is a cone. Clearly it contains with any vector
$(a,b)$ all vectors $\lambda(a,b)$, where $\lambda\geq0$. Also it is convex,
since the function $\psi_{t}(a,b)$ is concave, which follows from the
representation
\[
\psi_{t}(a,b)=\sum_{i=1}^{m}[(1-\lambda_{t,i}^{+})(R_{t,i}a^{i}-b^{i}%
)]-\sum_{i=1}^{m}[(\lambda_{t,i}^{-}+\lambda_{t,i}^{+})(R_{t,i}a^{i}%
-b^{i})_{-}],
\]
where the first sum is a linear function of $(a,b)$ and the second sum is a
convex function of $(a,b)$.

The model of a financial market we deal with corresponds to the von
Neumann-Gale dynamical system with the cones $X_{t}(\omega)$ specified by
(\ref{Xt-example}) and the cones $Z_{t}(\omega)$ given by \ref{Example2-Gt1}.
Paths in this dynamical system are self-financing trading strategies. Rapid
paths generalize benchmark strategies \cite{Platen2006,PlatenHeath2006} and
numeraire portfolios \cite{Long1990}.

We provide conditions that guarantee that the present model satisfies
conditions (\textbf{A1})-(\textbf{A4}), and so Theorem 1 is valid for it.
Define $\Lambda_{t,i}^{+}(\omega)=1-\lambda_{t,i}^{+}$ and $\Lambda_{t,i}%
^{-}(\omega)=1+\lambda_{t,i}^{-}$. Let us introduce the following conditions.

(\textbf{B1}) For each $t$, there exist constants $\underline{R}_{t}$,
$\overline{R}_{t}$, $\underline{\Lambda}_{t}$, $\overline{\Lambda}_{t}$ such
that $0<\underline{R}_{t}\leq R_{t,i}(\omega)\leq\overline{R}_{t}$,
$0<\underline{\Lambda}_{t}\leq\Lambda_{t,i}^{+}(\omega)$, $\Lambda_{t,i}%
^{-}(\omega)\leq\overline{\Lambda}_{t}$ for all $i$, $\omega$.

(\textbf{B2}) For each $t$, we have $\mu_{t}>\nu_{t}$ where
\[
\nu_{t}:=\max\{(\overline{\Lambda}_{t+1}\overline{R}_{t+1})/(\underline
{\Lambda}_{t+1}\underline{R}_{t+1});\overline{\Lambda}_{t}/\underline{\Lambda
}_{t}\}.
\]

\textbf{Proposition 9}. \textit{Let conditions (\textbf{B1}) and (\textbf{B2})
hold. Then the cones }$X_{t}(\omega)$\textit{\ satisfy condition (\textbf{A1})
and the cones }$Z_{t}(\omega)$\textit{\ satisfy conditions (\textbf{A2}%
)-(\textbf{A4}).}

To prove Proposition 9 we need the following auxiliary result.

\textbf{Lemma 1.} \textit{Let conditions (\textbf{B1}) and (\textbf{B2}) hold.
Then }

\textit{(a) For each }$t$\textit{\ there exists a constant }$C_{t}^{1}%
>0$\textit{\ such that for every }$a\in X_{t}(\omega)$\textit{\ the inequality
}$|a_{+}|-\nu_{t}|a_{-}|\geq C_{t}^{1}|a|$\textit{\ holds. }

\textit{(b) For each }$t$\textit{\ there exists a constant }$C_{t}^{2}%
$\textit{\ such that if }$a\in X_{t-1}(\omega)$\textit{, }$b\in X_{t}(\omega
)$\textit{\ and }$|b|\leq C_{t}^{2}|a|$\textit{, then }$(a,b)\in Z_{t}%
(\omega)$\textit{. }

\textit{Proof}\textbf{.} (a) By virtue of (\ref{Xt-example}), we have
$X_{t}(\omega)\subseteq\tilde{X}_{t}=\{a\in\mathbb{R}^{m}:\mu_{t}|a_{-}%
|\leq|a_{+}|\}$, where $\tilde{X}_{t}\cap(-\tilde{X}_{t})=\{0\}$ since
$\mu_{t}>1$. Observe that the continuous function $h_{t}(a)=|a_{+}|-\nu
_{t}|a_{-}|$ is strictly positive on the compact set $\hat{X}_{t}:=\tilde
{X}_{t}\cap\{a:|a|=1\}$. Indeed, since $h_{t}(a)\geq(\mu_{t}-\nu_{t})|a_{-}|$
on $\tilde{X}_{t}$, the equality $h_{t}(a)=0$ would imply $|a_{-}|=0$, and
hence $|a_{+}|=h_{t}(a)=0$, so that $|a|=0$. Then $h_{t}(a)$ attains a
strictly positive minimum on $\hat{X}_{t}$, which can be taken as $C_{t}^{1}$.

(b) Let $b\in X_{t}(\omega)$. It is straightforward to check that for any
numbers $x,y$ we have $(x-y)_{+}\geq x_{+}-y_{+}$ and $(x-y)_{-}\leq
x_{-}+y_{+}$. Using this, for any $a\in X_{t-1}(\omega)$ we obtain%

\[
\psi_{t}(a,b)\geq\sum_{i}(\Lambda_{t,i}^{+}R_{t,i}a_{+}^{i}-\Lambda_{t,i}%
^{-}R_{t,i}a_{-}^{i})-\sum_{i}(\Lambda_{t,i}^{+}+\Lambda_{t,i}^{-})b_{+}^{i}%
\]%
\[
\geq\underline{\Lambda}_{t}\underline{R}_{t}|a_{+}|-\overline{\Lambda}%
_{t}\overline{R}_{t}|a_{-}|-2\overline{\Lambda}_{t}|b_{+}|\geq\underline
{\Lambda}_{t}\underline{R}_{t}(|a_{+}|-\nu_{t-1}|a_{-}|)-2\overline{\Lambda
}_{t}|b|
\]%
\[
\geq C_{t-1}^{1}\underline{\Lambda}_{t}\underline{R}_{t}|a|-2\overline
{\Lambda}_{t}|b|.
\]
Assertion (b) will be valid for the constant $C_{t}^{2}:=C_{t-1}^{1}%
\underline{\Lambda}_{t}\underline{R}_{t}/(2\overline{\Lambda}_{t})$, since if
$|b|\leq C_{t}^{2}|a|$, then $\psi_{t}(a,b)\geq0$, implying $(a,b)\in Z_{t}$.

The proof is complete.

\textit{Proof of Proposition 9}. Let us check (\textbf{A1}). Consider the
non-random cone $\tilde{X}_{t}:=\{a\in\mathbb{R}^{m}:\mu_{t}|a_{-}|\leq
|a_{+}|\}$, so that $X_{t}(\omega)\subseteq\tilde{X}_{t}$. Put $q_{t}=e$,
where $e=(1,...,1)\in\mathbb{R}^{m}$. We can see that $q_{t}\in X_{t}^{\ast
}(\omega)$ since for any $a=(a^{1},...,a^{m})\in X_{t}(\omega)$, we have%
\[
q_{t}a=\sum_{i=1}^{m}a^{i}=|a_{+}|-|a_{-}|\geq(\mu_{t}-1)|a_{-}|\geq0.
\]
Observe that the continuous function $q_{t}a=\sum_{i=1}^{m}a^{i}$ is strictly
positive on the compact set $\hat{X}_{t}=\tilde{X}_{t}\cap\{a:|a|=1\}$.
Indeed, since $q_{t}a\geq(\mu_{t}-1)|a_{-}|$ on $\tilde{X}_{t}$, the equality
$q_{t}a=0$ would imply $|a|=0$. Then $q_{t}a$ attains a strictly positive
minimum $Q_{t}\leq1$ on $\hat{X}_{t}$. Define $H_{t}=Q_{t}^{-1}$. Hence, for
any $a\in X_{t}(\omega)$ we get%
\[
H_{t}^{-1}|a|\leq q_{t}a\leq H_{t}|a|,
\]
which implies that assumption (\textbf{A1}) is satisfied.

Condition (\textbf{A2}) follows from statement (b) of Lemma 1 since for any
$a\in X_{t-1}(\omega)$, $0\leq C_{t}^{2}|a|$ and so $(a,0)\in Z_{t}(\omega)$.

To prove (\textbf{A3}), let $(a,b)\in Z_{t}(\omega)$. Since for any numbers
$x,y$ we have $(x-y)_{+}\leq x_{+}+y_{-}$ and $(x-y)_{-}\geq y_{+}-x_{+}$, we obtain%

\[
0\leq\psi_{t}(a,b)\leq\sum_{i}(\Lambda_{t,i}^{+}+\Lambda_{t,i}^{-}%
)R_{t},ia_{+}^{i}+\sum_{i}(\Lambda_{t,i}^{+}b_{-}^{i}-\Lambda_{t,i}^{-}%
b_{+}^{i})
\]%
\begin{equation}
\leq2\overline{\Lambda}_{t}\overline{R}_{t}|a|+\overline{\Lambda}_{t}%
|b_{-}|-\underline{\Lambda}_{t}|b_{+}|\leq2\overline{\Lambda}_{t}\overline
{R}_{t}|a|-C_{t}^{1}\underline{\Lambda}_{t}|b|,\label{psi-estimate}%
\end{equation}
where in the last inequality, we used that $b\in X_{t}(\omega)$ and according
to statement (a) of Lemma 1, we have $\underline{\Lambda}_{t}|b_{+}%
|-\overline{\Lambda}_{t}|b_{-}|\geq\underline{\Lambda}_{t}(|b_{+}|-\nu
_{t}|b_{-}|)\geq C_{t}^{1}\underline{\Lambda}_{t}|b|$. This implies the
validity of (\textbf{A3}) with the constant $K_{t}=2\overline{\Lambda}%
_{t}\overline{R}_{t}/(C_{t}^{1}\underline{\Lambda}_{t})$.

Now we will prove condition (\textbf{A4}). Let $\mathring{x}=(1,...,1)\in
\mathbb{R}^{m}$. Put $\mathring{z}_{t}=(\mathring{x},\mathring{y}_{t})$ with
$\mathring{y}_{t}=(C_{t}^{2}/2)\mathring{x}$. Observe that there exists
$\delta_{t}>0$ such that $\mathbb{B}(\mathring{z}_{t},\delta_{t}%
)\subset\mathbb{R}_{+}^{2m}$ and therefore $\mathbb{B}(\mathring{z}_{t}%
,\delta_{t})\subset X_{t-1}\times X_{t}$. Since $|\mathring{y}_{t}|<C_{t}%
^{2}|\mathring{x}|$, statement (b) of Lemma 1 implies $\mathring{z}_{t}\in
Z_{t}$. Hence, $\mathring{z}_{t}\ $and $\delta_{t}$ satisfy condition
(\textbf{A4}).

The proof is complete.

\end{document}